\documentclass{svjour3}
\usepackage{comment}
\usepackage{enumerate}
\usepackage{bm}
\usepackage{amsmath,amssymb}

\newcommand{\R}{{\mathbb{R}}}

\newcommand{\cc}{{\bm{c}}}

\renewcommand{\ll}{{\bm{l}}}

\newcommand{\kk}{{\bm{k}}}

\newcommand{\xx}{{\bm{x}}}
\newcommand{\yy}{{\bm{y}}}
\newcommand{\xh}{{\widehat{\bm{x}}}}
\newcommand{\yh}{{\widehat{\bm{y}}}}

\newcommand{\CY}{{C_Y}}
\newcommand{\CE}{{C_E}}

\newcommand{\ih}{{ \widehat{i}^{\lambda_0} }}
\newcommand{\kh}{{ \widehat{k}^{\lambda_0} }}

\newcommand{\inp}[1]{{\langle{#1}\rangle}}

\newcommand{\Ih}{{ \widehat{I} }}
\newcommand{\q}{{ |q| }}
\newcommand{\sign}{{ \operatorname{sign} }}

\begin{document}

\title{An Accurate and Quadrature-Free Evaluation of Multipole Expansion of Functions Represented by Multiwavelets}
\titlerunning{Multipole Expansion of Multiwavelet Representations}

\author{Jae-Seok Huh} 

\institute{
	Jae-Seok Huh \at
	Computational Mathematics Group, Computer Science and Mathematics Division, Oak Ridge National Laboratory, Oak Ridge, TN 37831 \\
	Tel.: 1-865-574-3133\\
	Fax: 1-865-241-4811 \\
	\email{huhj@ornl.gov}
	} 

\date{Received: date / Accepted: date}

\maketitle

\begin{abstract}
We present formulae for accurate numerical conversion between functions represented by multiwavelets and their multipole/local expansions
with respect to the kernel of the form, $e^{-\lambda r}/r$ (cf.~\cite{GH02}).
The conversion is essential for the application of fast multipole methods for functions represented by multiwavelets.
The corresponding separated kernels exhibit near-singular behaviors at large $\lambda$.
Moreover, a multiwavelet basis function oscillates more wildly as its degree increases.
These characteristics in combination render any brute-force approach based on numerical quadratures impractical.
Our approach utilizes the series expansions of the modified spherical Bessel functions
and the Cartesian expansions of solid harmonics so that the multipole--multiwavelet conversion matrix can be evaluated like a special function.
The result is a quadrature-free, fast, reliable, and machine precision accurate scheme to compute the conversion matrix
with predictable sparsity patterns.
\end{abstract}


\section{Introduction}
Multiwavelets~\cite{A90,ABGV02} developed originally by  Alpert generalize Haar wavelets with piecewise polynomial
scale functions up to any given degree, hence, enjoy higher order of accuracy in the representation of sufficiently
smooth functions. 
A detailed discussion on the sparse representation of differential operators and exponential operators
for evolution equations can be found in~\cite{ABGV02}.
Recent advances in multi-resolution algorithms for integral operators can be found in various articles including
\cite{BM04,HFYGB04,YFGHB04,YFGHB04b},
which are based on the dimensional kernel separation via representation of kernels by weighted
sum of Gaussians, thus, applicable to a variety of kernels in arbitrary dimensions.

In this paper, we focus ourselves on a more traditional, but very popular and well-studied fast convolution algorithm
-- the fast multipole method. Our goal is to establish the connection between multiwavelet representation and
the fast multipole method. The problem can be reduced to finding the multipole expansion of multiwavelet basis functions.

Following one of the most recent version of the fast multipole algorithm presented in~\cite{GH02},
we begin with the kernel,
\begin{equation}
	G(\xx,\yy) = \frac{e^{-\lambda \| \xx-\yy \|}}{\| \xx-\yy \|}
\label{eqn:the_kernel}
\end{equation}
where $\lambda$ is a non-negative real number. This kernel is the fundamental solution of the linear differential operator,
\begin{equation}
	\Delta-\lambda^2 ,
\end{equation}
which appears in various applications involving damped Coulomb forces.
The resulting potential is also known as Yukawa potential in nuclear physics.
In the derivation of the scheme, we rely on the multipole expansion formula for strictly positive $\lambda$.
However, it turns out that the first term of the series representation of multipole expansion (with $\lambda>0$)
corresponds to the case of $\lambda=0$.
Hence, readers can assume that our scheme can be applied to \emph{the} Laplacian kernel also.

\section{The Multipole Expansion of Multiwavelet Basis Functions}

\subsection{The Multipole Expansion}
Denote by $\inp{\cdot,\cdot}$ the $L^2$ inner product of complex functions on a bounded domain in $\R^3$.
Let $\phi$ be a scalar \emph{source} function supported on the domain.
The potential generated by $\phi$ (outside of its support) is given by the following multipole series,
\begin{align}
	\Phi(\xx)
	&= \inp{G(\xx,\yy),\phi(\yy)} \\
	&= \sum_{p=0}^\infty \sum_{q=-p}^p M_p^q \, k_p(\lambda\|\xx\|) \, Y_p^q(\xh)
\end{align}
where $\xh=\xx/\|\xx\|$ is identified to a point on $S^2$.
The multipole coefficients $M_p^q$ are given by
\begin{equation}
	M_p^q = 8 \, \lambda \, \inp{i_p(\lambda\|\yy\|)Y_p^q(\yh),\phi(\yy)} .
\label{eqn:multipole}
\end{equation}
The functions $i_p$ and $k_p$ are the modified spherical Bessel and Hankel functions,
\begin{align}
	i_p(r) &= \sqrt{\frac{\pi}{2r}} \, I_{p+1/2}(r) \\
	k_p(r) &= \sqrt{\frac{\pi}{2r}} \, K_{p+1/2}(r) .
\end{align}
Since $i_p(r)$ and $k_p(r)$ exhibit exponential growth and decay respectively,
an algorithm based on the above formula is likely to experience an overflow/underflow.
To avoid the issue, as suggested in~\cite{GH02}, we replace
$i_p$ and $k_p$ in the formula by their scaled forms,
\begin{align}
	\ih_p (\lambda\,r) &= i(\lambda\,r)/\lambda_0^p \\
	\kh_p (\lambda\,r) &= k(\lambda\,r) \cdot \lambda_0^p.
\end{align}
Assuming $O(r)=1$ by appropriate geometric scaling, the appropriate choice of $\lambda_0$ is $\lambda$ itself, 
however, in order to maintain the generality, we keep clear notational distinction between them.

\begin{remark}
To avoid any confusion, the normalization of spherical harmonics $Y_p^q$ used in this paper needs to be clearly stated
before we begin any formulation.
We utilize exactly the same form presented in~\cite{GH02}, that is,
\begin{equation}
	Y_p^q (\yh) = \sqrt{\frac{2p+1}{4\pi}} \sqrt{\frac{(p-|q|)!}{(p+|q|)!}} \, P_p^{|q|}(\cos\theta) \, e^{i\,q\,\phi} .
\end{equation}
An obvious advantage of using this form is that $Y_p^{-q} = \overline{Y_p^q}$, hence, in (\ref{eqn:multipole}),
$Y_p^q$ appears without minus sign in front of $q$.
We can observe, later in this paper, that this property provides us with a better symmetry/sparsity pattern
of the multipole expansion matrix.
\end{remark}

\subsection{Symmetries}
The key equation in the above formula is the multipole expansion~(\ref{eqn:multipole}).
For notational simplicity, in this paper, we omit the constant $8\,\lambda$, which can be multiplied afterward.
We denote the product of $\ih_p$ and $Y_p^q$ in~(\ref{eqn:multipole}) by $Q_p^q$;
\begin{equation}
\label{eqn:Qpq}
	Q_p^q(\lambda,\yy) \equiv \ih_p(\lambda\|\yy\|) \, Y_p^q(\yh)
\end{equation}
The function $Q_p^q$ should be replaced with the regular solid harmonics when $\lambda=0$.
In the following section, we can observe that the first term of the series representation of $Q_p^q$
is the regular solid harmonics.
It is obvious that
\begin{equation}
	Q_p^q(\lambda,\alpha\,\yy) = Q_p^q(\alpha\,\lambda,\yy).
\end{equation}
The function $Q_p^q$ enjoys two useful symmetries:
firstly, from the normalization of $Y_p^q$ employed in this paper, it follows that
\begin{equation}
\label{eqn:Qpq:symm:1}
	Q_p^{-q} = \overline{Q_p^q} .
\end{equation}
Secondly, by change of variables $\phi\rightarrow\pi/2-\phi$ in $Y_p^q$, we can obtain
\begin{equation}
\label{eqn:Qpq:symm:2}
	Q_p^q(\lambda,y_2,y_1,y_3) = i^q \overline{Q_p^q}(\lambda,y_1,y_2,y_3) .
\end{equation}
As a result, a multipole expansion matrix presented in this paper possesses similar symmetries,
which we utilize to reduce the number of elements we have to compute.

\subsection{Multiwavelet Basis Functions}
In this section, we briefly introduce multiwavelet representation of functions.
A detailed discussion on the subject can be found in \cite{ABGV02}.
Denote by $\kk$ non-negative multi-indices and
by $\phi^\kk$ multi-dimensional orthonormal polynomials of degree $k_i$ in $i$th dimension.
We further assume that the generating functions $\phi^\kk$ are constructed
by the Cartesian product of 1-d orthonormal polynomials on $[-1,1]$, 
\begin{equation}
\label{eqn:e_def}
	\phi^\kk(\yy) = \sum_{i=1}^d \, \phi^{k_i}(y_i).
\end{equation} 

The above  $\phi^\kk$ generate the orthonormal multiwavelet basis functions at arbitrary level $n=0,1,\ldots$,
and translation characterized by multi-indices $\ll=(l_1,\ldots,\l_d)$ with $l_i=0,\ldots,2^n-1$
by formula,
\begin{equation}
	\phi^\kk_{n,\ll}(\xx) =
	\begin{cases}
	\sqrt{2}^{d(n+1)} \phi^\kk(2\,(2^n \xx - \ll) - 1) & \text{on $b_{(n,\ll)}$} \\
	0 & \text{elsewhere}
	\end{cases}
\end{equation}
where $b_{(n,\ll)}=\prod_{i=1}^d [2^{-n}l_i,2^{-n}(l_i+1)]$.
In this paper, we take $[0,1]^d$ ($= b_0$ by definition) as the computational domain.
Beware that we assume that the 1-d generating functions $\phi^{k_i}$ are orthogonal polynomials
defined on $[-1,1]$ (not on $[0,1]$).
This choice of \emph{unshifted} orthogonal polynomials as the generating functions
is to simplify the notations in multipole related formulae; we have to evaluate multipole expansions
with respect to the center of each $b_{(n,\ll)}$.

\begin{remark}
The term \emph{orthogonal polynomials} can be a source of confusion, which we need to clarify
before we present any related formula.
By the term, we mean a sequence of polynomials $\phi^k$ of degree $k$ orthogonal
to each other with respect to an underlying weighting function
(as in ``orthogonal polynomials and quadratures").
Since a non-trivial weighting function loses its meaning under scaling, readers may think
$\phi^k$ a synonym of (normalized) Legendre polynomial of degree $k$.
This limitation of generating functions to orthogonal polynomials greatly simplifies
the resulting formulae and makes the conversion matrix more sparse.
\end{remark}

\begin{remark}
There can be different choices of polynomial basis which are mutually orthogonal
such as the \emph{interpolating basis} presented in~\cite{ABGV02}.
Conversion between them and Legendre-generated basis is not complicated.
The advantage (by symmetry and sparsity) of using orthogonal polynomials exceeds
the additional cost of basis conversion.
\end{remark}

\subsection{The Multipole Expansion of $\phi^\kk_{(n,\ll)}$}
For any $p=0,1,\ldots$ and $q=-p,\ldots,p$, define $E^{(p,q)}_\kk(n,\lambda)$ by
\begin{equation}
	E^{(p,q)}_\kk(n,\lambda) = \inp{Q_p^q(\lambda,\yy-\cc_{(n,\ll)}),\phi^\kk_{(n,\ll)}(\yy)}_{b_{(n,\ll)}}
\end{equation}
where $\cc_{(n,\ll)}$ is the center of $b_{(n,\ll)}$.
The above equation gives the multipole coefficient $M_p^q$ (without $8\,\lambda$)
of (\ref{eqn:multipole}) with respect to the center, $\cc_{(n,\ll)}$.
The inner product can be scaled and translated to the standard domain $[-1,1]^3$,
\begin{align}
	E^{(p,q)}_\kk(n,\lambda)
	&= \frac{1}{\sqrt{2}^{3(n+1)}} \, \inp{Q_p^q(\lambda,2^{-(n+1)}\yy),\phi^\kk(\yy)} \\
	&= \frac{1}{\sqrt{2}^{3n}} \, E^{(p,q)}_\kk(0,\lambda/2^n)
\end{align}
where
\begin{equation}
	E^{(p,q)}_\kk(0,\lambda_n) = \frac{1}{\sqrt{2}^3} \, \inp{Q_p^q(\lambda_n/2,\cdot),\phi^\kk} .
\label{eqn:E0}
\end{equation}
Thus, we are required to evaluate (\ref{eqn:E0}) for arbitrary $\lambda_n$ which depends on
$\lambda$ and the level, $n$.
Viewing $(p,q)$ as a row multi-index and $\kk$ as a column multi-index,
$E^{(p,q)}_\kk$ acts as the conversion (multipole expansion) matrix for multiwavelet represented functions;
let $s^\kk_{(n,\ll)}$ be multiwavelet coefficients for a fixed $(n,\ll)$.
The multipole expansion centered at $\cc_{(n,\ll)}$ of the function
$\sum_\kk s^\kk_{(n,l)}\,\phi^\kk_{(n,\ll)}$ is given by the matrix-vector multiplication,
\begin{equation*}
	\sum_\kk E^{(p,q)}_\kk(n,\lambda) \, s^\kk_{(n,\ll)} .
\end{equation*}

The same matrix can be used for the conversion from a local expansion to its multiwavelet representation.
Consider a local expansion with the coefficients $L_p^q$,
\begin{equation}
	\Phi(\xx) = \sum_{p=0}^\infty \sum_{q=-p}^p \, L_p^q \, Q_p^q(\lambda,\xx) .
\end{equation}
The projection of $\Phi$ on to the span of the multiwavelet basis is, from the orthonormality, given by
\begin{equation}
	s^\kk_{(n,\ll)} = \inp{\overline{\Phi},\phi^\kk_{(n,\ll)}}
	= \sum_{(p,q)} \, \overline{E^{(p,q)}_\kk}(n,\lambda) \, L_p^q ,
\end{equation}
i.e., by the multiplication with the conjugate transpose of $E^{(p,q)}_\kk$.

\subsection{Symmetries}
Recall the symmetries of $Q_p^q$.
The following two conditions are the immediate consequences of (\ref{eqn:Qpq:symm:1}) and (\ref{eqn:Qpq:symm:2}).
\begin{equation}
\label{eqn:Epqk:symm:1}
	E^{(p,-q)}_\kk(n,\lambda) = \overline{E^{(p,q)}_\kk}(n,\lambda)
\end{equation}
and
\begin{equation}
\label{eqn:Epqk:symm:2}
	E^{(p,q)}_{(k_2,k_1,k_3)}(n,\lambda) = (-i)^q \, \overline{E^{(p,q)}_{(k_1,k_2,k_3)}}(n,\lambda) .
\end{equation}

In a later section, we will show that, depending on $(p,q,\kk)$, (1) $E^{(p,q)}_\kk$ is either real or pure imaginary
and (2) has pre-determined sparsity patterns.
Combined with the above symmetries, we recommend the following storage for the multipole expansion
matrices. For each level $n$, we compute $E^{(p,q)}_\kk$ for $q\ge0$ and $k_2>k_1$, and store non-negative $q$ portion
of the matrices in two sparse matrices,
one for real and the other for imaginary.
Separated storage is simple and advantageous in the implementation;
(i) since each element is either real or imaginary, the two sparse matrices have disjoint index sets;
it does not require any additional storage or computation cost due to duplicated indices.
(ii) For rows with $q<0$, the multiplication can be omitted; for a complex vector $s^\kk_{(n,\ll)}$,
\begin{equation*}
	\sum_{\kk} E^{(p,-q)}_{\kk,\Re/\Im} \, s^\kk_{(n,\ll)} =
	\pm \sum_{\kk} E^{(p,q)}_{\kk,\Re/\Im} \, s^\kk_{(n,\ll)}
\end{equation*}
with the negative sign for imaginary matrix.

\subsection{Numerical Issues}
There are three major numerical issues which make the evaluation of $E^{(p,q)}_\kk(n,\lambda)$ non-trivial.

\begin{enumerate}[(1)]
\item {\bf Non-homogeneity of $\bm{Q_p^q(\lambda,\cdot)}$}:
Unlike regular solid harmonics, $Q_p^q$ are not homogeneous.
Since we cannot extract $\lambda$ out of the integral,
we have to build different $E^{(p,q)}_\kk(n,\lambda)$ depending on $\lambda$ and $n$,
which rules out the possibility of utilizing a precomputed table.
Since they are not even polynomials, there is no simple quadrature
which produces the exact integral.
Any naive approach using adaptive quadrature becomes impractical for the following reasons.

\item {\bf Rapid growth of $\bm{Q_p^q(\lambda,\cdot)}$}:
The function $i_p(\lambda\,r)$ grows exponentially.
Scaling by using $\ih_p(\lambda\,r)$ helps preventing overflow.
However, the function still exhibits near singularity for large $\lambda$.

\item {\bf Oscillating behavior of $\bm{\phi^\kk}$}:
Although they are polynomials, $\phi^k$ have all their zeros on $(-1,1)$.
Hence, for large $k$, an adaptive integrator will encounter with highly oscillating integrands.
\end{enumerate}

From the above characteristics, any adaptive integration requires a large amount of computation,
or simply it fails to converge especially for highly oscillating cases.
Beware that, to build $N$ conversion matrices (up to depth level $N-1$) for the multipole expansion
(up to degree $P$) of functions represented by multiwavelets (of degree up to $K$),
we have to compute $O(N\times P^2\times K^3)$ elements!

Our approach begins with rewriting the $Q_p^q$ in a series form.
Each term in the resulting power series involves a regular solid harmonics weighted by
an even power of $\|\xx\|$, hence, a homogeneous polynomial in $\R^3$.
The Cartesian expansion of this polynomial can be explicitly written and its projection
on multiwavelet basis can be obtained exactly without using any numerical quadrature.

\section{The Series Form}

In this section, we present a series representation of the multipole expansion matrix, $E^{(p,q)}_\kk(n,\lambda)$.
We begin with the identity,
\begin{equation}
	I_\alpha(r) = \sum_{m=0}^\infty \frac{1}{m! \, \Gamma(m+\alpha+1)} \left( \frac{r}{2} \right)^{2m+\alpha} .
\end{equation}
Utilizing $\Gamma(m+1/2)=\sqrt{\pi}\,(2m)!/(4^m\,m!)$,
we obtain
\begin{equation}
	\ih_p(\lambda\,r)
	= r^p \left(\frac{\lambda}{\lambda_0}\right)^p
	\sum_{m=0}^\infty \frac{1}{m! \, (2m+2p+1)!!} \, \left(\frac{\lambda\,r}{\sqrt{2}}\right)^{2m} .
\end{equation}
Therefore $Q_p^q$ can be written in power series given by
\begin{equation}
	Q_p^q(\lambda,\xx)
	= \left(\frac{\lambda}{\lambda_0}\right)^p
	\sum_{m=0}^\infty \frac{1}{m! \, (2m+2p+1)!!} \, \left(\frac{\lambda^2}{2}\right)^m
	\|\xx\|^{2m+p}\,Y_p^q(\xh)
	.
\end{equation}

Define $R_{p,m}^q$ by
\begin{equation}
	R_{p,m}^q(\xx) = \CY(p,q)^{-1} \, \|\xx\|^{2m+p}\,Y_p^q(\xh)
\end{equation}
where
\begin{equation}
	\CY(p,q) = \sqrt{\frac{2p+1}{4\pi}}\,\sqrt{\frac{(p-|q|)!}{(p+|q|)!}} ,
\end{equation}
Notice that the function $R_{p,m}^q$ is a regular solid harmonics multiplied by $\|\xx\|^{2m}$,
hence, a homogeneous polynomial of degree $(2m+p)$.
The factor $\CY(p,q)^{-1}$ simplifies the Cartesian expansion of $R_{p,m}^q$, which we
introduced in the following section.

From the above series representation, $E^{(p,q)}_\kk(0,\lambda_n)$ is given by
\begin{align}
	E^{(p,q)}_\kk(0,\lambda_n)
	&= \frac{1}{\sqrt{2}^3} \inp{Q_p^q(\lambda_n/2,\cdot),\phi^\kk} \nonumber \\
	&= \CE(\lambda_n/\lambda_0,p,q)
	\sum_{m=0}^\infty A_m(p) \, \left(\frac{\lambda_n^2}{8}\right)^m
	I_m(p,q,\kk)
\label{eqn:series}
\end{align}
where
\begin{equation}
	\CE(\lambda_n/\lambda_0,p,q)
	= \frac{\CY(p,q)}{\sqrt{2}^3\,(2p+1)!!} \left(\frac{\lambda_n}{\lambda_0}\right)^p ,
\end{equation}

\begin{equation}
	A_m(p) = \frac{(2p+1)!!}{m! \, (2m+2p+1)!!} ,
\end{equation}
\begin{equation}
	I_m(p,q,\kk) = \int_{[-1,1]^3} \overline{R_{p,m}^q}(\xx) \, \phi^\kk(\xx) \, d\xx .
\end{equation}
The factor $(2p+1)!!$ in $A_m(p)$ is added for the normalization, $A_0(p)=1$.

In (\ref{eqn:series}), $\lambda_n$ is now taken out of the integral.
We will observe that the $\lambda$-independent term, $I_m(p,q,\kk)$, can be further reduced
to a finite sum of products of 1-d integrals with two integer parameters, namely, $\Ih_k^l$.
We can construct $\Ih_k^l$ exactly without using any numerical quadrature
via the recurrence relation of orthogonal polynomials.
Our strategy is to tabulate $\Ih_k^l$ and use the table to evaluate the series (\ref{eqn:series})
for various $\lambda$, $n$, $p$, $q$, and $\kk$.

\subsection{Properties of $I_m$ and the Convergence Criterion}

Most of the properties of $I_m$ presented in this section will be explained in detail in~\S\ref{sec:Im}.
For a more comprehensive presentation, we think it would be more appropriate to
discuss the behavior of (\ref{eqn:series}) prior to the presentation of detailed formulae for $I_m$.
Followings are the summary of the relevant properties:
\begin{enumerate}
\item $I_m(p,q,\kk)$ is either real or pure imaginary depending only on $\kk$. 
\item $I_m(p,q,\kk)=0$ if $2m<k_x+k_y+k_z-p$. 
\item Sign of $I_m(p,q,\kk)$ is determined by $q$ only and is independent of $m$.
\item $|I_m(p,q,\kk)|$ is monotonically increasing as $m$ increases.
\item $I_{m+1}(p,q,\kk)/I_m(p,q,\kk)\rightarrow3$ as $m\rightarrow\infty$.
\end{enumerate}

\begin{remark}
Property (5) can be supported by the following estimate:
Since $|\CY^{-1}\overline{Y_p^q}|=|P_p^{|q|}(\cos\theta)\,e^{i\,q\,\phi}|\le1$,
\begin{align*}
|I_m(p,q,\kk)| 
&\le \|\phi^\kk\|_\infty \, \int_{[-1,1]^3} \|\xx\|^{2m+p} d\xx \\
&< 4\pi \, \|\phi^\kk\|_\infty \, \int_0^{\sqrt{3}} r^{2m+p+2} dr = \frac{12\pi\,\sqrt{3}^p\,3^m}{(2m+p+3)} \, \|\phi^\kk\|_\infty
\end{align*}
\end{remark}

Thus, the series consists of two parts: $A_m$ decreasing factorially
and $(\lambda^2/8)^m\,I_m$ which behaves asymptotically $\sim (3\lambda^2/8)^m$.
Their product $C_m=A_m (\lambda^2/8)^m I_m$ has a fixed sign for a fixed $(p,q,\kk)$
independently of $m$.
Hence, the partial sum of the series increases (decreases) monotonically
to the upper (lower) bound which is potentially huge in the absolute sense.
The non-alternating feature of the series 
suppresses any necessity of considerations of cancellation errors, and
suggest the following simple convergence criterion:
for given absolute and relative tolerances $\epsilon_a$ and $\epsilon_r$, stop the summation if
 \begin{equation}
 \label{eqn:conv}
 	|C_M| < \epsilon_a
	\quad\text{or}\quad
	|C_M| < \epsilon_r \, \left| \sum_{m=0}^M C_m \right|
	.
 \end{equation}
 
We can numerically observe that the number of terms to convergence $M$ is $O(\lambda)$
in a conservative estimation.
For example, for $\epsilon_r=10^{-16}$, $M\sim \lambda$ and slightly smaller if $\lambda$ is large; e.g.,
when $\lambda\sim300$, $M\sim200$.
The condition (2) combined with the convergence criterion
provides us with additional sparsity of $E^{(p,q)}_\kk$;
if $2M<k_x+k_y+k_z-p$, the corresponding $E^{(p,q)}_\kk$ can be considered to be zero.

\begin{remark}
We represent $E^{(p,q)}_\kk$ like a special function of $\lambda$ with exponential growth.
The number of terms $M$ can grow indefinitely as $\lambda$ increases.
Although, in many practical applications, $\lambda$ are quite limited and 
$\lambda_n=\lambda/2^n$ decreases as the depth of the multiwavelet representation increases,
a more complete algorithm requires an asymptotic expansion of $E^{(p,q)}_\kk$ with respect to $\lambda$.
Yet, we haven't found a closed formula for the asymptotic expansion, which is an on-going work.
\end{remark}

\section{The Formula for  $I_m$ and the Sparsity Pattern}
\label{sec:Im}

In this section, we present an explicit Cartesian expansion form of $\overline{R_{p,m}^q}$ in $I_m$.
Each term can be written as a product of 1-d integrals which can be evaluated
exactly by the recurrence relations of the orthogonal polynomials $\phi^k$.
We begin with the series form of the spherical harmonics.
With the Rodrigues' formula, the associated Legendre functions $P_p^{|q|}$ in $Y_p^q$
can be written as
\begin{align}
	P_p^\q (z)
	&= \frac{(-1)^\q}{2^p \, p!} \, (1-z^2)^{\q/2} \, \frac{d^{p+\q}}{dz^{p+\q}} \, (z^2-1)^p \nonumber \\
	&= \frac{(-1)^\q}{2^p} \, (1-z^2)^{\q/2} \, \sum_{\nu=0}^{\lfloor\frac{p-\q}{2}\rfloor}
		\, \frac{(-1)^\nu(2p-2\nu)!}{\nu! \, (p-\nu)! \, (p-\q-2\nu)!} \, z^{(p-\q-2\nu)}
	.
\end{align}
Hence, using notations $\xx=(x,y,z)$, $r=\|\xx\|$, and $s = \sign(q)$,
\begin{multline}
	\overline{R_{p,m}^q}(\xx)
	= r^{2m+p} P_p^\q(z/r) \left(\frac{x-s\,i\,y}{\sqrt{r^2-z^2}}\right)^\q \\
	= \frac{(-1)^\q}{2^p} \, (x-s\,i\,y)^\q
		\sum_{\nu=0}^{\lfloor\frac{p-\q}{2}\rfloor}
		\frac{(-1)^\nu(2p-2\nu)!}{\nu! \, (p-\nu)! \, (p-\q-2\nu)!} \, r^{2(m+\nu)} \, z^{(p-\q-2\nu)}
\end{multline}
By expanding $(x-s\,i\,y)^\q$ and $r^{2(m+\nu)}$, we obtain
\begin{multline}
	= (s\,i)^\q \\
	\sum_{\mu=0}^\q (s\,i)^\mu a_\mu
	\sum_{\nu=0}^{\lfloor\frac{p-\q}{2}\rfloor} b_\nu
	\sum_{\alpha=0}^{m+\nu} c_{\nu\alpha} \, z^{(p-\q+2m-2\alpha)}
	\sum_{\beta=0}^\alpha d_{\alpha\beta} \, y^{(\q+2\beta-\mu)} \, x^{(2\alpha-2\beta+\mu)}
\label{eqn:expansion}
\end{multline}
where the coefficients are given by
\begin{align}
	a_\mu &= \begin{pmatrix} q \\ \mu \end{pmatrix} \\
	b_\nu &= \frac{(-1)^\nu}{2^\nu} \frac{(2p-2\nu-1)!!}{\nu! \, (p-q-2\nu)!} \\
	c_{\nu\alpha} &= \begin{pmatrix} m+\nu \\ \alpha \end{pmatrix} \\
	d_{\alpha\beta} &= \begin{pmatrix} \alpha \\ \beta \end{pmatrix}
\end{align}
where we use the definition, $(-1)!!=0!!=1$.

\subsection{The Formula}
\label{sec:Im:formula}

From (\ref{eqn:expansion}), we obtain our final formula for $I_m$:
\begin{equation}
\label{eqn:Im}
	I_m(p,q,\kk) = (s\,i)^\q \, I_m^{(1)}(p,q,\kk) + (s\,i)^{\q+1} \, I_m^{(2)}(p,q,\kk)
\end{equation}
where
\begin{multline}
\label{eqn:Im1}
	I_m^{(1)}(p,q,\kk) =
	\sum_{\mu=0}^{\lfloor\frac{\q}{2}\rfloor} (-1)^\mu a_{2\mu}
	\sum_{\nu=0}^{\lfloor\frac{p-\q}{2}\rfloor} b_\nu
	\\ \cdot
	\sum_{\alpha=0}^{m+\nu} c_{\nu\alpha} \, \Ih_{k_z}^{(p-\q+2m-2\alpha)}
	\sum_{\beta=0}^\alpha d_{\alpha\beta} \, \Ih_{k_y}^{(\q+2\beta-2\mu)} \, \Ih_{k_x}^{(2\alpha-2\beta+2\mu)}
	,
\end{multline}
\begin{multline}
\label{eqn:Im2}
	I_m^{(2)}(p,q,\kk) =
	\sum_{\mu=0}^{\lfloor\frac{\q-1}{2}\rfloor} (-1)^\mu a_{2\mu+1}
	\sum_{\nu=0}^{\lfloor\frac{p-\q}{2}\rfloor} b_\nu
	\\ \cdot
	\sum_{\alpha=0}^{m+\nu} c_{\nu\alpha} \, \Ih_{k_z}^{(p-\q+2m-2\alpha)}
	\sum_{\beta=0}^\alpha d_{\alpha\beta} \, \Ih_{k_y}^{(\q+2\beta-2\mu-1)} \, \Ih_{k_x}^{(2\alpha-2\beta+2\mu+1)}
	,
\end{multline}
and
\begin{equation}
	\Ih_k^l = \int_{-1}^1 \zeta^l \phi^k(\zeta) d\zeta .
\end{equation}

Note that $I_m^{(1)}$ and $I_m^{(2)}$ are real functions and,
with factors $(s\,i)^\q$ and $(s\,i)^{\q+1}$ respectively,
they determine real and imaginary parts of $I_m$ separately.
The above representation of $I_m$ by two separate parts $I_m^{(1)}$ and $I_m^{(2)}$
is to signify the following very useful fact: at least, one of $I_m^{(1)}$ and $I_m^{(2)}$ vanishes
for any $(p,q,\kk)$ independently of $m$,
which implies that
\begin{equation*}
	\text{$E^{(p,q)}_\kk$ is either a real or a pure imaginary.}
\end{equation*}
Moreover, depending on the parameter $(p,q,\kk)$, many of $I_m$ vanish,
which results in the nice sparsity of the multipole expansion matrix.
These properties are the immediate consequence of the following properties of orthogonal polynomials.

\paragraph{{\rm(1)} {\bf Oddity}}
Any orthogonal polynomial $\phi^k$ with symmetric domain and weight is even (odd)
if the degree $k$ is even (odd). Hence,
\begin{equation*}
	\Ih_k^l = 0 \quad\text{if $(l+k)$ is odd.}
\end{equation*}
Since $m$ always appears in the equation with the factor of 2,
any consequence of the oddity condition is $m$-independent;
the resulting sparsity of $E^{(p,q)}_\kk$ is pre-determined by $(p,q,\kk)$ only
(independently of the level $n$ and $\lambda$).
We can observe that
\begin{enumerate}[(a)]
\item $I_m^{(2)}=0$ if $q=0$.
\item $I_m^{(1)}=0$ if $k_x$ is odd or $(\q+k_y)$ is odd or $(p+\q+k_z)$ is odd.
\item $I_m^{(2)}=0$ if $k_x$ is even or $(\q+k_y)$ is even or $(p+\q+k_z)$ is odd.
\end{enumerate}
Therefore,
\begin{equation}
	E^{(p,q)}_\kk = 0
	\quad\text{if}\quad
	\begin{cases}
	\text{$(k_z+p+\q)$ is odd} \\
	\text{$(k_x+k_y+\q)$ is odd} \\
	\text{$q=0$ and at least one of $k_x$ and $k_y$ is odd}
	\end{cases}
\end{equation}
Suppose $E^{(p,q)}_\kk\ne0$ from the above test.
Then, the oddity of $k_x$ must be the same as the oddity of $(\q+k_y)$, which results in
\begin{equation}
	I_m = c \cdot \begin{cases}
	I_m^{(2)} & \text{if $k_x$ is odd} \\
	I_m^{(1)} & \text{if $k_x$ is even}
	\end{cases}
	\quad
	c = \begin{cases}
	(-1)^{\lfloor \frac{\q}{2} \rfloor} \,\sign(q)\, i & \text{if $k_y$ is odd} \\
	(-1)^{\lceil \frac{\q}{2} \rceil} & \text{if $k_y$ is even}
	\end{cases}
\end{equation}
Notice, $E^{(p,q)}_\kk$ is real (imaginary) if $k_y$ is even (odd).
	
The following table illustrates the sparsity of the multipole expansion matrix for parameters:
$0 \le p \le 10$, $0 \le q \le p$, and $0 \le k_{x,y,z} \le 10$.
We can observe that about a quarter of elements are non-zeroes. (See Table~\ref{tab1}.)

\begin{table}[hbt]
\begin{center}
\begin{tabular}{ccc} \hline
	total elements & real non-zeroes & imaginary non-zeroes \\ \hline
	87846 & 12186 (13.9\%) & 8450 (9.6\%) \\ \hline
\end{tabular}
\end{center}
\caption{$\lambda$-independent sparsity estimated by the oddity condition.}
\label{tab1}
\end{table}

\paragraph{{\rm(2)} {\bf Moment condition}}
Recall the moment conditions satisfied by orthogonal polynomials.
\begin{equation*}
	\Ih_k^l = 0 \quad\text{if $l < k$.}
\end{equation*}
Consider a term with a fixed set of indices $(p,q,\kk,\mu,\nu,\alpha,\beta)$ in (\ref{eqn:Im}).
The term vanishes if
\begin{equation*}
	k_z > p - \q + 2m - 2\alpha \quad\text{or}\quad
	k_y > \q + 2\beta - \mu \quad\text{or}\quad
	k_x > 2\alpha - 2\beta + \mu
\end{equation*}
which is true if
\begin{equation*}
	k_x + k_y + k_z > p + 2m .
\end{equation*}
Thus,
\begin{equation}
	I_m(p,q,\kk) = 0 \quad\text{if}\quad 2m < k_x + k_y + k_z - p
\end{equation}
This condition is $m$-dependent, hence, cannot be used to pre-determine the sparsity pattern
of $E^{(p,q)}_\kk$. However, it still can affect the sparsity for a given $\lambda$;
suppose that the convergence criterion (\ref{eqn:conv}) is satisfied at $M$ for $2M<k_x+k_y+k_z-p$.
Then, the corresponding $E^{(p,q)}_\kk(n,\lambda)$ is effectively zero.
Since $M\sim\lambda$ and $E^{(p,q)}_\kk(n,\lambda)=constant\cdot E^{(p,q)}_\kk(0,\lambda/2^n)$,
the multipole expansion matrix becomes more sparse as $\lambda$ decreases and as $n$ increases.
Table~\ref{tab2} shows the number of vanishing elements (for $n=0$)
among those predicted to be non-zeroes by the oddity condition.
With parameters, $0 \le p \le 10$, $0 \le q \le p$, and $0 \le k_{x,y,z} \le 10$,
the number of elements is 87846 and the numbers of non-zero real and imaginary elements
(predicted by the oddity condition) are 12186 and 8450 respectively (same as the above example).
Tolerances are $\epsilon_a=\epsilon_r=10^{-16}$.

\begin{table}[hbt]
\begin{center}
\begin{tabular}{ccc} \hline
	$\lambda$ & additional real zeroes & additional imaginary zeroes \\ \hline
	1 & 9567 & 6679 \\ \hline
	2 & 8813 & 6154 \\ \hline
	4 & 7478 & 5235 \\ \hline
	6 & 6340 & 4439 \\ \hline
	8 & 5439 & 3775 \\ \hline
	10 & 4630 & 3203 \\ \hline
	50 & 1309 & 851 \\ \hline
	100 & 1301 & 848 \\ \hline
	200 & 1273 & 835 \\ \hline
	300 & 1251 & 824 \\ \hline
\end{tabular}
\end{center}
\caption{$\lambda$-dependent sparsity estimated from the moment condition.}
\label{tab2}
\end{table}

The additional sparsity decreases as $\lambda$, hence $M$, increases.
There are two factors which controls $M$ and, hence, the the additional sparsity
-- the absolute tolerance $\epsilon_a$ and the relative tolerance $\epsilon_r$.
Among them, the contribution of $\epsilon_r$ decreases rapidly
and becomes quite negligible when $\lambda\gg\max\,k_i$.
However, the contribution of $\epsilon_a$ is persistent.
From the table, we can observe that the additional sparsities by $\epsilon_a$
are $\sim1200$ for the real matrix and $\sim800$ for the imaginary matrix.

Also, the moment condition enhances the computational efficiency slightly.
We can simply skip the evaluation of $I_m$ if $m<\lceil(k_x+k_y+k_z-p)/2\rceil$.
Like the sparsity by $\epsilon_r$, the effect decreases quite rapidly as $\lambda$ increases.
However, in practical applications of fast multipole and multiwavelet representation, 
the levels of terminal boxes where we need to perform the expansion is likely to be high.
Hence, the additional sparsity by $\epsilon_r$ should not be considered insignificant.

\subsection{The Laplacian Kernel {\rm($\bm{\lambda=0}$)}}
In this case, $Q_p^q$ becomes simply the regular solid harmonics.
The corresponding $E^{(p,q)}_\kk(n,0)$ is just the first term ($m=0$) of the series form (\ref{eqn:series})
with an appropriate adjustment of the constant factor.
In this case, the moment condition becomes
\begin{equation*}
	E^{(p,q)}_\kk = 0 \quad\text{if}\quad p < k_x + k_y + k_z
\end{equation*}
which results in a more sparse multipole expansion matrices.
With the same condition as the previous examples, $\max p = \max k_i = 10$,
the number of non-zero elements are given in Table~\ref{tab3}.

\begin{table}[hbt]
\begin{center}
\begin{tabular}{ccc} \hline
	total elements & real non-zeroes & imaginary non-zeroes \\ \hline
	87846 & 1512 (1.72\%) & 1001 (1.14\%) \\ \hline
\end{tabular}
\end{center}
\caption{Sparsity of the case with $\lambda=0$.}
\label{tab3}
\end{table}

We can observe that the resulting matrices are very sparse --
only less than 3\% of total elements are non-zeroes.
This example illustrates the efficiency of the multipole expansion
on multiwavelet representations based on
\emph{orthogonal} (almost synonymously in this paper, Legendre) polynomials.

\subsection{Recurrence Relations for $\Ih_k^l$}
Finally, we present the algorithm to build the table of $\Ih_k^l$
required for the evaluation of $I_m$.
Let $p\le p_{\max}$, $k_i\le k_{\max}$, and $M\le M_{\max}$ (for a given $\lambda$).
Then, the required size of table is $(2M_{\max}+p_{\max})\times k_{\max}$, and
the half of the elements are zero by the oddity condition.

Each element $\Ih_k^l$ can be calculated by the identical recurrence relations to those of the
orthogonal polynomials $\phi^k$. Recall any sequence orthogonal polynomials satisfy a three term recurrence relation
of the form,
\begin{equation}
	\phi^{k+1} = (\alpha_k \, \zeta + \beta_k) \, \phi^k - \gamma_k \, \phi^{k-1} .
\end{equation}
It immediately follows that
\begin{equation}
	\Ih_{k+1}^l = \alpha_k \, \Ih_k^{l+1} + \beta_k \, \Ih_k^l - \gamma_k \, \Ih_{k-1}^l .
\end{equation}
From the oddity condition, $\Ih_{k+1}^l\ne0$ if and only if $\Ih_k^l=0$,
and $\Ih^l_k=0$ for $l<k$.
Hence, for $k < l$,
\begin{equation}
	\Ih_k^{l+1} = 
	\begin{cases}
		0 & \text{if $(l+k)$ is even} \\
		\alpha_k^{-1} \, \Ih_{k+1}^l + \alpha_k^{-1} \gamma_k \, \Ih_{k-1}^l & \text{if $(l+k)$ is odd}
	\end{cases}
\end{equation}
and
\begin{equation}
	\Ih_k^k = \alpha_k^{-1} \, \gamma_k \, \Ih_{k-1}^{k-1}
	\quad
	(\because \Ih_{k+1}^{k-1}=0).
\end{equation}
The recurrence relation can be evaluated from the initial data,
\begin{equation}
	\Ih_0^l = a_0 \, \frac{1+(-1)^l}{l+1} \quad\text{and}\quad
	\Ih_1^l = a_1 \, \frac{1-(-1)^l}{l+2}
\end{equation}
where $\phi^0=a_0$ and $\phi^1=a_1\zeta$.

For \emph{normalized} Legendre polynomials, the coefficients are given by
\begin{equation}
	\alpha_k^{-1} = \frac{k+1}{2k+1} \, \sqrt{\frac{2k+1}{2k+3}}
	\quad\text{and}\quad
	\alpha_k^{-1} \gamma_k = \frac{k}{2k+1}  \, \sqrt{\frac{2k+1}{2k-1}} .
\end{equation}

\section{Results and Conclusions}
Most of the implementation is done very faithfully with the formulae presented in this paper.
The only special treatment is that, in order to suppress the accumulation of round-off errors affecting the result,
we used higher precision floating point arithmetics for internal calculations including the table for $\Ih^l_k$;
for example, to generate matrices with 64bit double precision, we employed 80bit long(extended)-double arithmetics.
The computational cost is governed by the number of terms to be added, $m$, and is not significantly affected by
the augmented internal precision.
By comparing with values obtained by applying adaptive numerical integrator to (\ref{eqn:e_def}),
we could validate the presented formula.
When $\lambda$ or $\kk$ is only moderately large, an adaptive integrator typically fails to converge
since the integrand of (\ref{eqn:e_def}) becomes near-singular or highly oscillating.
Thus, the presented formula can be viewed as a reliable way to evaluate (\ref{eqn:e_def}) (or a similar form of integral)
when a typical numerical quadrature is not applicable due to the near-singularity and/or the oscillation of the integrand.
It is also observed that the computing time of building $\Ih_k^l$ table is negligible
compared to the computing time of $E^{(p,q)}_\kk$.
We summarize the contributions of this paper as follows.
\begin{enumerate}
\item We presented a method to build the multipole expansion matrices for functions represented by
multiwavelets.
\item The presented method does not involve any numerical quadrature and based entirely on
 a series representation like a special function of $\lambda$.
\item The proposed scheme generates highly accurate multipole conversion matrices
stably and reliably for a wide range of parameters $(p,q,\kk)$ and $\lambda$.
\end{enumerate}

\end{document}